\documentclass[12pt,a4paper]{article}
\usepackage{amsmath,amsthm,amsfonts,amssymb}

\newtheorem{theorem}{Theorem}[section]

\newtheorem{corollary}[theorem]{Corollary}

\newcommand{\Z}{\mathbb Z}
\newcommand{\Q}{\mathbb Q}
\renewcommand{\:}{\colon}
\renewcommand{\>}{\rightarrow}

\title{Criteria of measure-preserving for $p$-adic dynamical systems in terms of the van der Put basis}

\author{Andrei Khrennikov and Ekaterina Yurova\\ International Center for Mathematical Modelling\\
in Physics and Cognitive Sciences\\ 
Linnaeus University V\"axj\"o, S-35195, Sweden\\
andrei.khrennikov@lnu.se}

\begin{document}
\maketitle

\begin{abstract} This paper is devoted  to (discrete) $p$-adic dynamical systems,  an important domain of algebraic and arithmetic 
dynamics. We consider the following open problem from theory of $p$-adic dynamical systems. Given continuous function $f\:\Z_p\>\Z_p.$ Let us represent it via special convergent series, namely van der Put series. How can one specify whether this function is measure-preserving or not for an arbitrary $p?$
In this paper, for any prime $p,$ we present a complete description of all compatible measure-preserving functions in the additive form representation. 
In addition we prove the criterion in terms of coefficients with respect to the van der Put basis determining whether a compatible function 
$f\:\Z_p\>\Z_p$ preserves the Haar measure.
\end{abstract}

\section{Introduction}

Algebraic and arithmetic dynamics are actively developed fields of general theory of dynamical systems. The bibliography collected by Franco Vivaldi \cite{Vivaldi}
contains 216 articles and books; extended bibliography also can be found in books of Silverman \cite{Silverman} and Anashin and Khrennikov \cite{ANKH}. Theory of dynamical systems in fields of $p$-adic numbers and their algebraic extensions is  an important part of 
algebraic and arithmetic dynamics, see, e.g., \cite{Albeverio1}--\cite{Vivaldi2} (the complete list of reference would be very long; 
hence, we refer to \cite{Vivaldi}--
\cite{ANKH}, \cite{Khrennikov and Nilsson}). As in general theory of dynamical systems, problems of ergodicity and measure preserving play fundamental roles in
theory of $p$-adic dynamical systems, see \cite{ANKH}--\cite{Arrowsmith2},  \cite{Khrennikov and Nilsson}, \cite{Parry}. 
Traditionally studies in these domains of p-adic dynamics were restricted to analytic (mainly polynomial) or at least smooth
maps $f:{\bf Q}_p \to {\bf Q}_p,$ where ${\bf Q}_p$ is the field of $p$-adic numbers. However, the internal mathematical development of theory of $p$-adic 
dynamical systems as well as applications to cryptography \cite{ANKH} stimulated the interest to nonsmooth dynamical maps. An important class of (in general) 
nonsmooth maps is given by {\it Lipschitz one functions.} In cryptographic applications such functions are called {\it compatible}. We shall use this 
terminology in this paper. 

The main  mathematical tool used in this paper is the representation of the function by the {\it van der Put series} which is actively used in
$p$-adic analysis, see e.g. Mahler \cite{Mahler} and Schikhof \cite{Shihov}. 
Marius van der Put introduced this series in his dissertation ``Alg�bres de fonctions continues p-adiques'' at Utrecht Universiteit in 1967, \cite{vdp}.
There are numerous results in studies of functions with zero derivatives, antiderivation \cite{Shihov} obtained using van der Put series. Later van der Put basis was adapted to the case of $n$-times continuously differentiable functions in one and several variables \cite{DeSmedt}. First results on applications of the van der Put series
in theory of $p$-adic dynamical systems, the problems of ergodicity and measure preserving,  were obtained in \cite{DAN}.  
The present paper is the first attempt to use van der Put basis to examine such property as measure-preserving
of (discrete) dynamical systems in a space
of $p$-adic integers $\Z_p$ for an arbitrary prime $p.$ 

Note that the van der Put basis differs fundamentally from previously used ones (for
example, the monomial and Mahler basis \cite{Erg})  which are
related to the algebraic structure of $p$-adic fields. The
van der Put basis is related to the zero dimensional
topology of these fields (ultrametric structure), since it
consists of characteristic functions of $p$-adic balls. In other words,
the basic point in the construction of this basis is the
continuity of the characteristic function of a $p$-adic
ball.

In this paper, we present  a description of all compatible measure-preserving functions by using the  {\it additive form representation,}
Theorem 4.1. In the additive form criteria, a compatible measure-preserving function  is decomposed into a sum of two functions. The first one is an arbitrary compatible function - ``free'' part, and the second one is a compatible function of a special type. This ``special'' function is given by the van der Put basis, where the coefficients are defined via an arbitrary set of substitution on the set of nonzero residues modulo $p$ and one substitution modulo $p.$

The  additive form representation,
Theorem 4.1,  is based on the criterion of measure-preserving in terms of coefficients of the van der Put series, see
Theorem 2.1. It was announced in \cite{DAN}, \cite{Ya}. In this paper we give its proof.
As an example of its application,  we show how known results  on the description of certain classes of compatible measure-preserving functions  can be obtained from Theorem 3.2.
Namely, the classes of compatible $2$-adic functions and uniformly differentiable functions modulo $p.$

\section{Criterion of measure-preserving}

Let $p>1$ be an arbitrary prime number. The ring of $p$-adic integers is denoted by  the symbol $\Z_p.$ The $p$-adic valuation
is denoted by $\vert \cdot\vert_p.$ We remind that this valuation satisfies the {\it strong triangle inequality}:
$$
\vert x+y\vert_p \leq \max [\vert x\vert_p, \vert y\vert_p].
$$
This is the main distinguishing property of the $p$-adic valuation inducing essential departure from the real or complex analysis
(and hence essential difference of  $p$-adic dynamical systems from  real and complex dynamical systems).

We shall use the terminology of  papers \cite{DAN},  \cite{Ya}. 

Namely, van der Put series are defined in the following way. Let $f\:\Z_p\>\Z_p$ be a continuous function. Then there exists a unique sequence of $p$-adic coefficients
$B_0,B_1,B_2, \ldots$ such that 
\begin{eqnarray}
f(x)=\sum_{m=0}^{\infty}
B_m \chi(m,x) 
\end{eqnarray}
for all $x \in \Z_p.$
Here the \emph{characteristic function} $\chi(m,x)$ is given by
\begin{displaymath}
\chi(m,x)=\left\{ \begin{array}{cl}
1, \mbox{if}
\ \left|x-m \right|_p \leq p^{-n} \\
0,  \mbox{otherwise}
\end{array} \right.
\end{displaymath}
where $n=1$ if $m=0$, and $n$ is uniquely defined by the inequality 
$p^{n-1}\leq m \leq p^n-1$ otherwise (see Schikhof's book \cite{Shihov} for detailed presentation 
of theory of van der Put series). 

The \emph{van der Put coefficients} $B_m$ are
related to the values of $f$ as follows.
Let 
$m=m_0+ \ldots +m_{n-2} p^{n-2}+m_{n-1} p^{n-1}$ be the representation of $m$ in the $p$-ary number system, i.e., 
 $ m_j\in \left\{0,\ldots ,p-1\right\}$, $j=0,1,\ldots,n-1$ and $m_{n-1}\neq 0.$ Then

\begin{displaymath}
B_m= \left\{ \begin{array}{cl}
f(m)-f(m-m_{n-1} p^{n-1}), \mbox{if}\; \; 
m\geq p \\
f(m),  \mbox{otherwise}
\end{array} \right.
\end{displaymath}

 Let $f\:\Z_p\>\Z_p$ be a function  and
let $f$ satisfy the {\it Lipschitz condition with constant $1$} (with respect to the $p$-adic valuation $\left|\cdot \right|_p):$
 $$\left| f(x)- f(y)\right|_p \leq \left|x-y \right|_p$$
for all $x,y \in \Z_p.$ 

We state again that a mapping of an algebraic system $A$ to itself is called {\it compatible} if it preserves all the congruences of $A.$ 
It is easy to check that a map $f: \Z_p \to \Z_p$ is Lipschitz one iff it is compatible (with respect to $\rm{mod}\; p^k, k=1,2,...$
congruences). 

The space $\Z_p$ is equipped by the natural probability measure, namely, the Haar measure $\mu_p$ normalized
so that $\mu_p(\Z_p) = 1.$ 

Recall that a mapping $f\:\mathbb S\rightarrow \mathbb S$ of a measurable space $\mathbb S$ with a probability measure $\mu$ is called measure preserving if 
$\mu(f^{-1}(S))=\mu(S)$ for each measurable subset $S\subset \mathbb S.$

We say that a compatible function $f\:\Z_p\>\Z_p$ is bijective modulo $p^k$ if the induced mapping $x\mapsto f(x)\bmod p^k$ is a permutation on $\Z/p^k\Z.$ 
It was shown in [9] (see also [1, Section 4.4]) that a compatible function $f\:\Z_p\>\Z_p$ is measure-preserving if and only if it is bijective modulo $p^k$ for all $k=1,2,3,\ldots .$

\begin{theorem}
\label{v1}
Let $f\: \Z_p \> \Z_p$ be a compatible function and 
\[
f(x)=\sum_{m=0}^\infty p^{\left\lfloor \log_p m\right\rfloor} b_m \chi (m,x)
\]
 be the van der Put representation of this function, where $b_m \in \Z_p, m=0,1,2,\ldots .$ Then $f(x)$ preserves the Haar's measure iff
\begin{enumerate}
\item $b_0,b_1,\ldots, b_{p-1}$ establish a complete set of residues modulo $p,$ i.e. the function $f(x)$ is bijective modulo $p;$
\item 
\[
b_{m+p^k}, b_{m+2p^k},\ldots ,b_{m+(p-1)p^k}
\]
for any $m=0,\ldots, p^k-1$  are all nonzero residues modulo $p$ for $k=2,3,\ldots .$
\end{enumerate}
\end{theorem}

Proof.
By the induction for $k=1, 2,\ldots$ we show that the function $f$ is bijective modulo $p^k.$ For $k=1$ it is true by the first condition of the theorem, i.e. 
$b_i\equiv f(i) \bmod p$ for $i=0,1,\ldots,p-1.$
Assume that $f$ is bijective modulo $p^k.$ Let us show that the function $f$ is bijective modulo $p^{k+1}.$ In other words, we should show that the comparison 
$f(x)\equiv \acute{t} + p^kt \bmod p^{k+1}$ has a unique solution for any $\acute{t}\in \left\{0,\ldots, p^k-1\right\}$ and $t=0,1,\ldots,p-1.$
By the induction hypothesis the comparison $f(x)\equiv \acute{t} \bmod p^k$ has a unique solution $\acute{x}\in \left\{0,\ldots, p^k-1\right\}.$
Then to check bijective property of the function $f \bmod p^{k+1}$ it is enough to show that for a given value $\acute{t}\in \left\{0,\ldots, p^k-1\right\}$ the 
comparison 
\begin{eqnarray}
\label{eq:sr1}
f(\acute{x} +p^k x)\equiv \acute{t}+p^k t \bmod p^{k+1}
\end{eqnarray}
has unique solution with respect to $x\in \left\{0,\ldots, p-1\right\}$ for any $t\in \left\{0,\ldots, p-1\right\}.$

For every $\acute{x}\in \left\{0,\ldots, p^k-1\right\}$ we set a function
\[
\varphi _{\acute{x}}(h)=
\left\{
\begin{array}{cl}
b_{\acute{x}+p^k h} \bmod p, & \mbox{$h\neq 0$} \\
0, & \mbox{$h=0$}
\end{array}
\right.
\]

which is defined and valued in the residue ring  modulo $p.$

To calculate values of the function $f$ using the van der Put representation we write a comparison (\ref{eq:sr1}) as
\begin{eqnarray}
\label{eq:sr2}
f(\acute{x} +p^k x)= f(\acute{x})+p^k \varphi_{\acute{x}}(x).
\end{eqnarray}

We take into account that $\acute{x}$ is unique solution of the comparison $f(x)\equiv \acute{t} \bmod p^k$ 
and assume $f(\acute{x})\equiv \acute{t}+p^k \xi \bmod p.$
Thus we transform the comparison (\ref{eq:sr2}) as 
\begin{eqnarray}
\label{eq:sr3}
\varphi _{\acute{x}}(x)\equiv t-\xi \bmod p.
\end{eqnarray}

Under the second condition of the theorem \ref{v1} the function $\varphi _{\acute{x}}$ is bijective on $\left\{0,\ldots,p-1\right\}.$
Then for any $t=0,\ldots,p-1$ the comparison (\ref{eq:sr3}) has unique solution on $\left\{0,\ldots,p-1\right\}.$
That is (\ref{eq:sr1}) has unique solution, and therefore the function $f$ is bijective modulo $p^k$ for any $k=1,2,\ldots .$
Thereby the function $f$ preserves measure by the Theorem 1.1 in \cite{Unif}.

Now let us prove the theorem in the opposite direction. Let the function $f$ preserves measure. By \cite{ANKH} 
$f$ is bijective modulo $p^{k+1}$ for any $k=1,2,\ldots .$ The first condition follows immediately from this result (here $k=0$).
The following comparisons have unique solution $\left(\acute{x};x\right),$ where $\acute{x}\in \left\{0,\ldots, p^k-1\right\}$ and $x\in \left\{0,\ldots, p-1\right\}$ 
for any $\acute{t}\in \left\{0,\ldots, p^k-1\right\}$ and $t\in \left\{0,\ldots, p-1\right\}:$
\begin{eqnarray}
f(\acute{x} +p^k x)\equiv \acute{t}+p^k t \bmod p^{k+1} \\
f(\acute{x})\equiv \acute{t} \bmod p^k.
\label{eq:sr4}
\end{eqnarray}

After transformations presented at the beginning of the proof we can see that the condition of uniqueness of the solution of comparisons (\ref{eq:sr4}) is equivalent to unique solvability of the comparison (\ref{eq:sr3}) with respect to $x\in \left\{0,\ldots,p-1\right\}$ for any $t\in \left\{0,\ldots,p-1\right\}.$
It means that the function $\varphi _{\acute{x}}$ is bijective on $\left\{0,\ldots,p-1\right\}.$
And therefore $b_{\acute{x}+p^k}, b_{\acute{x}+2p^k},\ldots ,b_{\acute{x}+(p-1)p^k}$
coinside with the set of all nonzero residues modulo $p.$

\bigskip

The formulation and proof of measure-preservation of the locally compatible $p$-adic functions are similar to the previous reasoning. 
Remind that locally compatible functions are ones satisfying the $p$-adic Lipschitz condition with
a constant of $1$ locally, i.e., in a suitable neighborhood
of each point from $\Z_p$, see \cite{DAN}.

\begin{corollary}
\label{v2}
Let $f\: \Z_p \> \Z_p$ be a locally compatible function and 
\[
f(x)=\sum_{m=0}^\infty p^{\left\lfloor \log_p m\right\rfloor} b_m \chi (m,x)
\]
 be the van der Put representation of this function, where $b_m \in \Z_p, m\geq N .$ Then $f(x)$ preserves the Haar's measure iff
\begin{enumerate}
\item the function $f(x)$ is bijective modulo $p^N;$
\item 
\[
b_{m+p^k}, b_{m+2p^k},\ldots ,b_{m+(p-1)p^k}
\]
for any $m=0,\ldots, p^k-1$  are all nonzero residues modulo $p$ for $k>N.$
\end{enumerate}
\end{corollary}

\section{Connection to known results}

Here we show how to use theorems above by proving some known results on description of compatible measure-preserving $p$-adic functions.

A compatible measure-preserving $2$-adic functions represented via the van der Put series have been described in papers, see for example \cite{DAN}, \cite{Ya},  and state that

\begin{theorem}
The function $f: \Z_2\rightarrow \Z_2$ is compatible and preserves  
the measure $\mu_p$ if and only if
it can be represented as
\[
f(x)= b_0\chi(0,x)+ b_1\chi(1,x)+  
\sum_{m=2}^{\infty}
2^{\left\lfloor \log_2 m \right\rfloor} b_m \chi(m,x),
\]
where $b_m\in \Z_2$ for $m=0,1,2,\ldots$, and
\begin{enumerate}
\item $b_0+b_1\equiv 1 \bmod 2$,
\item $\left| b_m \right|_2=1$, if $m\ge 2$.
\end{enumerate}
\end{theorem}

Turns out that this result immediately follows from Theorem \ref{v1}. Indeed, from the second condition of this theorem follows that $b_m\equiv 1 \bmod 2$ for $m\geq 2$ or, in another words, $\left| b_m \right|_2=1.$ The first condition of Theorem \ref{v1} means that $b_0+b_1\equiv 1 \bmod 2,$ which is equivalent to the first condition in the theorem mentioned above.

In papers \cite{ANKH}, \cite{Unif} outlined characterization of measure-preserving, uniformly differentiable modulo $p$ compatible $p$-adic functions $f: \Z_p\rightarrow \Z_p.$
We remind, see also Defenitions 3.27 and 3.28 from \cite{ANKH}, that the function $f: \Z_p\rightarrow \Z_p$ is uniformly differentiable modulo $p$ if for any $u\in \Z_p$ there exits 
a positive rational integer $N$ and $f^{'}_1(u)\in \Q_p$ such that for every $k\geq N$ and $h\in \Z_p$ we have $f(u+p^k h)\equiv f(u)+p^k h\cdot f^{'}_1 (u) \bmod p^{k+1}.$
Note that an uniformly differentiable modulo $p$ function is locally compatible as soon as $\left|f^{'}_1 (u)\right|_p\leq 1$ ($n=m=1$), according to the Proposition 3.41 \cite{ANKH}.

Therefore, we can see that the Theorem 4.45 \cite{ANKH} follows from our Corollary \ref{v2}. Indeed, from the definition of the uniformly differentiable function modulo $p$ and periodically function $f^{'}_1 (u)$ with period $p^N$, see Proposition 3.32 from \cite{ANKH}, follows that for sufficiently large $k$ there exists the comparison 
$b_{u+p^k h}\equiv p^k h\cdot f^{'}_1 (u) \bmod p^{k+1}.$ 
The set of the coefficients $b_{u+p^k h} \bmod p$ for $h=1,2,\ldots, p-1$ coincides with the set of all nonzero residues modulo $p$ iff $f^{'}_1 (u)\neq 0 \bmod p$ for any $u\in \Z_p.$ Thereby here the function $f$ preserves measure iff the function $f$ is bijective modulo $p^k$ for some $k\geq N$ and $f^{'}_1 (u)\neq 0 \bmod p,$ as well as stated Theorem 4.45 in \cite{ANKH}.

\section{Additive representation of compatible, measure-preserving $p$-adic functions}

By the Lemma 4.41 \cite{ANKH} we know that for arbitrary compatible function $g: \Z_p\rightarrow \Z_p$ the function $f(x)=d+cx+pg(x),$ $d,c\in \Z_p,$ $c\neq0 \bmod p$ preserves measure. So we see that  measure-preserving functions are linear combinations of some ``fixed'' measure-preserving function and ``arbitrary'' compatible function. Moreover, in the case $p=2$  necessary and sufficient conditions for measure-presurving property were obtained, \cite{ANKH}. It turns out that we can find a similar representation for all compatible measure-preserving functions for any prime $p$. Here as a ``fixed'' part the special class of measure-preserving functions appears.

\begin{theorem}
\label{v3}
Let $h\: \Z_p \> \Z_p$ is an arbitrary compatible function.
A compatible function $f\: \Z_p \> \Z_p$ preserves measure iff it can be represented as 
\[
f(x)=\xi (x)+ p\cdot h(x),
\]
where the functions $\xi (x)$ represented via the van der Put series is such that
\begin{eqnarray}
  \xi (x)= \sum_{m=0}^{p-1} G(i) \chi (i,x) + \sum_{k=1}^\infty \sum_{m=0}^{p^k-1} \sum_{i=1}^{p-1} g_m (i)p^k\cdot \chi (m+i\cdot p^k,x) = \nonumber\\
= \sum_{m=0}^{p-1} G(i) \chi (i,x) + \sum_{k=1}^\infty \sum_{m=0}^{p^k-1} \sum_{i=1}^{p-1} i\cdot p^k\cdot \chi (m+g_m^{-1} (i)\cdot p^k,x),
\end{eqnarray}
where $g_m$ is a substitution on the set $\left\{1,\ldots,p-1\right\}$ and $G$ is a a substitution on the set $\left\{0,1,\ldots,p-1\right\}.$
\end{theorem}

Proof.
By the hypothesis of this theorem the function $f$ is compatible. Let us represent $f$ via the following van der Put series. Here the van der Put coefficients are $B_m = p^{\left\lfloor \log_p m \right\rfloor} (b_m + p\cdot \tilde{b_m}),$ 
where $b_m,  \tilde{b_m} \in \Z_p.$  Then
\[
f(x)=\sum_{m=0}^\infty p^{\left\lfloor \log_p m \right\rfloor} b_m\cdot \chi (m,x) + p\cdot \sum_{m=0}^\infty p^{\left\lfloor \log_p m \right\rfloor} \tilde{b_m} \cdot \chi (m,x).
\]
By the theorem about compatibility  (i.e. the Lipschitz property with constant 1) 
of functions represented via the van der Put series, see \cite{Shihov},  \cite{DAN}, \cite{Ya},  the function $$h(x)=\sum_{m=0}^\infty p^{\left\lfloor \log_p m \right\rfloor} \tilde{b_m} \cdot \chi (m,x)$$ is compatible. Now we set $$\xi (x)= \sum_{m=0}^\infty p^{\left\lfloor \log_p m \right\rfloor} b_m\cdot \chi (m,x).$$
By Theorem \ref{v1} a compatible function $\xi (x)$ preserves measure iff it is bijective modulo $p,$ i.e. $b_i=G(i)$ for $i=0,1,\ldots, p-1,$ where 
$$b_{m+p^k}, b_{m+2p^k},\ldots ,b_{m+(p-1)p^k}$$ are all nonzero residues modulo $p$ for $k=1,2,3,\ldots $ and $m=0,\ldots, p^k-1.$
Let $g_m$ be a substitution on the set $1,2,\ldots,p-1$ such that $g_m(i)=b_{m+i\cdot p^k} \bmod p$ for $k=1,2,3,\ldots $ and $m=0,\ldots, p^k-1.$ 
Then 
\begin{eqnarray}
  \xi (x)= \sum_{m=0}^{p-1} G(i) \chi (i,x) + \sum_{k=1}^\infty \sum_{m=0}^{p^k-1} \sum_{i=1}^{p-1} g_m (i)p^k\cdot \chi (m+i\cdot p^k,x) = \nonumber\\
= \sum_{m=0}^{p-1} G(i) \chi (i,x) + \sum_{k=1}^\infty \sum_{m=0}^{p^k-1} \sum_{i=1}^{p-1} i\cdot p^k\cdot \chi (m+g_m^{-1} (i)\cdot p^k,x),
\end{eqnarray}
where $g_m^{-1}$ is an inverse substitution to $g_m.$

\subsection{Example of measure-preserving function in additive representation}
Now we consider an example of the compatible measure-preserving $p$-adic function constructed by using the additive representation.
For $k=0,1,2,\ldots$ we choose substitutions $g_m$ such that $g_0=g_1=\ldots =g_{p^k-1}=h_k.$
Denote via $\delta_k (x)$ the value of $k$-th $p$-adic digit in a canonical expansion of the number $x.$
Then  
\begin{eqnarray}
\sum_{m=0}^{p^k-1} \sum_{i=1}^{p-1} g_m (i)p^k\cdot \chi (m+i\cdot p^k,x) = \sum_{i=1}^{p-1} h_k(i) p^k\cdot \sum_{m=0}^{p^k-1} \chi (m+i\cdot p^k,x) = \nonumber\\
= \sum_{i=1}^{p-1} h_k(i) p^k\cdot I(\delta_k (x)=i),
\end{eqnarray}
where 
\[
I(\delta_k (x)=i) =
\left\{
\begin{array}{cl}
1,& \mbox{$\delta_k (x)=i$}\\
0,& \mbox{$\delta_k (x)\neq i$}.
\end{array}
\right.
\]

Then let us choose integer $s$ such that $GCD(s,p-1)=1.$ Substitutions $h_k$ defined on the set $\left\{1,2,\ldots,p-1 \right\}$ determine by the comparison $h_k(i)=i^s \bmod p$ for $k=1,2,\ldots .$ Using the equality $\sum_{i=0}^{p-1} G(i) \chi (i,x)=G(\delta_0(x))$ we represent the function $\xi(x)$ as
\[
\xi(x)= G(\delta_0(x)) + \sum_{k=1}^\infty p^k \cdot (\delta_k (x))^s \bmod p
\]
or
\[
\xi(x_0+x_1p+ \ldots + x_kp^k+\ldots)= G(x_0) + \sum_{k=1}^\infty p^k \cdot x_k ^s \bmod p.
\]

The substitution $G$ we define, let us say, $G(x_0)=p-1-x_0,$ where $x_0=\left\{0,\ldots, p-1\right\}.$
Then as a function $h(x)$ from the Theorem \ref{v3} we take the pseudo-constant function $h(x)=\sum_{k=0}^\infty x_k p^{2k}$ with $x=\sum_{k=0}^\infty x_k p^k \in \Z_p,$  example 26.4, p.74 \cite{Shihov}.
Finally we get that the function
\begin{eqnarray}
  f(x_0+x_1p+ \ldots + x_kp^k+\ldots)= (p-1)(1+x_0)+ \sum_{k=1}^\infty p^k \cdot x_k ^s \bmod p + \nonumber\\
+ \sum_{k=1}^\infty p^{2k+1} \cdot x_k = \nonumber\\
= (p-1)(1+x_0)+ \sum_{k=1}^\infty (x_{2k+1}^s \bmod p +x_k)\cdot p^{2k+1} + \nonumber\\
+ \sum_{k=1}^\infty p^{2k} \cdot x_{2k}^s \bmod p
\end{eqnarray}
preserves measure.

\section*{Acknowledgments}
 
This work is supported by the joint grant of Swedish and South-African Research Councils, ``Non-Archimedean analysis: from fundamentals 
to applications.''

\end{document}